\documentclass{amsart}

\usepackage{amsmath,amsthm}
\usepackage{amsfonts,amssymb}

\newtheorem{thm}{Theorem}[section]
\newtheorem{cor}[thm]{Corollary}
\newtheorem{lem}[thm]{Lemma}
\newtheorem{prop}[thm]{Proposition}

\theoremstyle{remark}

\def\a{{\alpha}}
\def\b{{\beta}}

 \def\CH{{\mathcal H}}
 \def\CJ{{\mathcal J}}
 \def\CV{{\mathcal V}}
 \def\RR{{\mathbb R}}

 \def\proj{\operatorname{proj}}
 \def\sspan{\operatorname{span}}
 \def\vol{\operatorname{vol}}

\newcommand{\wh}{\widehat}

\begin{document}

\title
{A family of Sobolev Orthogonal Polynomials on the Unit Ball}

\author{Yuan Xu}
\address{Department of Mathematics\\ University of Oregon\\
    Eugene, Oregon 97403-1222.}\email{yuan@math.uoregon.edu}

\date{\today}
\keywords{Sobolev orthogonal polynomials, several variables, unit ball}
\subjclass{42A38, 42B08, 42B15}
\thanks{The author has been supported by NSF Grant DMS-0201669}

\begin{abstract}
A family of orthonormal polynomials on the unit ball $B^d$ of $\RR^d$ 
with respect to the inner product 
$$
 \langle f,g \rangle
    = \int_{B^d}\Delta[(1-\|x\|^2) f(x)] \Delta[(1-\|x\|) g(x)] dx, 
$$
where $\Delta$ is the Laplace operator, is constructed explicitly. 
\end{abstract}

\maketitle

\section{Introduction}\label{Introduction}
\setcounter{equation}{0}

In a recent study \cite{A} on the numerical solution of the nonlinear Poisson 
equation $-\Delta u = f(\cdot, u)$ on the unit disk with zero boundary conditions, 
Atkinson \cite{A} asked the question of finding an explicit orthogonal basis for 
the inner product defined by
$$
  \langle f,g \rangle_\Delta := \frac{1}{\pi}\int_{B^2} 
        \Delta[(1-x^2-y^2) f(x,y)]\Delta[(1-x^2-y^2) g(x,y)] dx dy
$$
on the unit disk $B^2$ of the Euclidean plane, where $\Delta$ is the usual 
Laplace operator. The purpose of this note is to provide an answer for this 
question. 

We shall consider more generally the analogous inner product on the unit ball
$B^d$ in $\RR^d$. We call orthogonal polynomials with respect to such an 
inner product Sobolev orthogonal polynomials. In the theory of orthogonal 
polynomials of one variable, the name Sobolev is associated with polynomials 
that are 
orthogonal with respect to an inner product defined using both functions and 
their derivatives; see, for example, \cite{G} and the references therein. As
far as we know, Sobolev orthogonal polynomials have not been studied in 
the case of several variables.

Our main result, given in Section 2, is a family of orthonormal polynomials with 
respect to $\langle \cdot,\cdot \rangle_\Delta$ on $B^d$ that are constructed 
using spherical harmonics and Jacobi polynomials in Section 2. For $d =1$, 
orthogonal polynomials with respect to this inner product has been studied 
recently in \cite{H}. The explicit formula can be used to study further 
properties of the orthogonal basis. In particular, it turns out that the
orthogonal expansion of a function $f$ in this basis can be computed without 
involving the derivatives of $f$. This will be discussed in Section 3.

\section{Sobolev Orthogonal Polynomials}\label{Sobolev}
\setcounter{equation}{0}

For $x \in \RR^d$, let $\|x\|$ denote the usual Euclidean norm of $x$. 
The unit ball in $\RR^d$ is $B^d:=\{x: \|x\| \le 1\}$. Its surface is 
$S^{d-1}:=\{x:  \|x\| =1\}$. The volume of $B^d$ and the surface area of 
$S^{d-1}$ are denoted by $\vol(B^d)$ and $\omega_{d-1}$, respectively, 
$$
\vol(B^d) = \omega_{d-1}/d \quad\hbox{and}\quad 
    \omega_{d-1} = 2 \pi^{d/2}/\Gamma(d/2).
$$
Let $\Pi^d = \RR[x_1,\ldots,x_d]$ be the ring of polynomials in $d$ variables
and let $\Pi_n^d$ denote the subspace of polynomials of total degree at most 
$n$. We consider the inner product defined on the polynomial space by
$$
  \langle f,g \rangle_\Delta 
   : = \frac{1}{4 d^2 \vol(B^d)}
     \int_{B^d}\Delta[(1-\|x\|^2) f(x)] \Delta[(1-\|x\|^2) g(x)] dx. 
$$
The constants are chosen so that $\langle 1, 1 \rangle_\Delta =1$. As pointed 
out in \cite{A}, the inner product is well defined and positive definite on 
$\Pi^d$. Let $\CV_n^d(\Delta)$ denote the space of orthogonal 
polynomials with respect to $ \langle f,g \rangle_\Delta$. It follows from
the general theory of orthogonal polynomials in several variables (\cite{DX}) 
that the dimension of $\CV_n^d(\Delta)$ is $\binom{n+d-1}{d-1}$. A polynomial 
$P$ is in $\CV_n^d(\Delta)$ if it is orthogonal to all polynomials of lower 
degrees with respect to $\langle f,g \rangle_\Delta$. If $\{P_\alpha \}$ is a 
basis of $\CV_n^d(\Delta)$ and $\langle P_\alpha, P_\beta \rangle_\Delta = 0$
whenever $\alpha \ne \beta$, it is called a mutually orthogonal basis. If, 
in addition, $P_\alpha$ is normalized so that $\langle P_\alpha,P_\alpha 
\rangle_\Delta =1$, the basis is called orthonormal. Our objective in this
section is to find an explicit orthonormal basis for $\CV_n^d(\Delta)$. 

The presence of the Laplace operator suggests that we make use of harmonic 
polynomials, which are homogeneous polynomials that satisfy the equation 
$\Delta P =0$. Let $\CH_n^d$ denote the space of harmonic polynomials of 
degree $n$. It is well known that 
$$
\dim \CH_n^d = \binom{n+d-1}{d-1}-\binom{n+d-3}{d-1}:= \sigma_n.
$$ 
The restriction of $Y \in \CH_n^d$ on $S^{d-1}$ are called spherical harmonics.
They are orthogonal on $S^{d-1}$. We will use the spherical polar coordinates
$x = r x'$ for $x \in \RR^d$, $r\ge 0$, and $x' \in S^{d-1}$. For $Y \in \CH_n^d$ 
we use the notation $Y(x)$ to denote the harmonic polynomials and use $Y(x')$ 
to denote the spherical harmonics. This agrees with $x = r x'$ since $Y$ is a 
homogeneous polynomial, $Y(x) = r^nY(x')$. Throughout this paper, we use 
the notation $\{Y_\nu^n: 1 \le \nu \le \sigma_n\}$ to denote an orthonormal basis 
for $\CH_n^d$, that is,  
\begin{equation} \label{eq:harmonics}
\frac{1}{\omega_{d-1}} \int_{S^{d-1}} Y_\mu^n(x') Y_\nu^m(x') d\omega(x')
    = \delta_{\mu,\nu} \delta_{n,m} , \qquad 1 \le \mu,\nu\le \sigma_n, 
\end{equation}
where $d\omega$ stands for the surface measure on $S^{d-1}$. 
In terms of the 
spherical polar coordinates, $x = r x'$, $r > 0$ and $x' \in S^{d-1}$, the
Laplace operator can be written as
\begin{equation}\label{eq:Delta}
  \Delta = \frac{\partial^2}{\partial r^2}+ \frac{d-1}{r} 
      \frac{\partial}{\partial r} + \frac{1}{r^2} \Delta_0,
\end{equation}
where $\Delta_0$ is the spherical Laplacian on $S^{d-1}$. It is 
well-known that 
\begin{equation}\label{eq:Eigen}
\Delta_0 Y (x') = -n(n+d-2) Y(x'), \qquad Y \in \CH_n^d, \quad x'\in S^{d-1}.
\end{equation}

The spherical harmonics have been used to construct orthogonal polynomials on
the unit ball. For later use, let us mention an orthogonal basis with respect 
to the inner product 
$$
\langle f,g \rangle_\mu : = c_\mu \int_{B^d} f(x) g(x) W_\mu(x)dx,\qquad
   W_\mu(x) = (1-\|x\|^2)^{\mu}, 
$$
where $\mu > -1$ and $c_\mu$ is the normalization constant of $W_\mu$. Let 
$\CV_n^d(W_\mu)$ denote the space of orthogonal polynomials of degree $n$. 
A mutually orthogonal basis for $\CV_n(W_\mu)$ is given by (\cite{DX})
\begin{equation}\label{eq:Wmu-basis}
 P_{j,\nu}^n(W_\mu; x) = P_j^{(\mu, n-2j+\frac{d-2}{2})}(2\|x\|^2 -1) 
   Y_{\nu}^{n-2j}(x), \quad 0 \le j \le n/2,
\end{equation}
where $P_j^{(\a,\b)}$ denotes the Jacobi polynomial of degree $j$, which 
is orthogonal with respect to $(1-x)^\alpha(1+x)^\beta$ on $[-1,1]$, and
$\{Y_\nu^{n-2j}: 1 \le j \le \sigma_{n-2j}\}$ is a basis for $\CH_{n-2j}^d$. 

In view of \eqref{eq:Wmu-basis} we will look for a basis with respect to 
$\langle f,g \rangle_\Delta$ in the form of 
\begin{equation}\label{eq:Q-basis}
Q_{j,\nu}^n(x) = q_j(2\|x\|^2-1) Y_\nu^{n-2j}(x), \qquad 0 \le 2j \le n, 
     \quad Y_\nu^{n-2j} \in \CH_{n-2j}^d,
\end{equation} 
where $q_j$ is a polynomial of degree $j$ in one variable.

\begin{lem} \label{lem:2.1}
Let $Q_{j,\nu}^n$ be defined as above. Then
$$
 \Delta\left[(1-\|x\|^2)Q_{j,\nu}^n(x) \right] =
    4 \left(\CJ_\beta q_j\right)(2r^2-1) Y_\nu^{n-2j}(x), 
$$   
where $\beta = n-2j + \frac{d-2}{2}$ and
$$
(\CJ_\beta q_j)(s) = (1-s^2) q_j''(s) + (\beta -1 - (\beta+3) s )q_j'(s) 
     - (\beta+1) q_j(s). 
$$
\end{lem}

\begin{proof}
Using spherical-polar coordinates, \eqref{eq:Delta} and \eqref{eq:Eigen} 
show that 
\begin{align*}
 & \Delta\left[(1-\|x\|^2)Q_{j,\nu}^n(x) \right] = 
   \Delta\left[ (1-r^2) q_j(2r^2-1)r^{n-2j} Y_\nu^{n-2j}(x')\right] \\
 &\qquad  = 4 r^{n-2j} \left[4 r^2 (1-r^2)q_j''(2r^2-1) + 2((\beta +1)-
   (\beta+3)r^2) q_j'(2r^2-1) \right.\\
 & \qquad\qquad\qquad \qquad\qquad \left. 
   - (\beta+1) q_j(2r^2-1)\right] Y_\nu^{n-2j}(x').
\end{align*}
Setting $s \mapsto 2r^2-1$ gives the stated result. 
\end{proof}

\begin{lem} \label{lem:2.2}
Let $p_k^\beta \in \Pi_k:=\Pi_k^1$ be orthogonal with respect to the inner 
product
$$
  (f,g)_\beta := \int_{-1}^1 (\CJ_\beta f)(s)(\CJ_\beta g)(s)
    (1+s)^\beta ds, \qquad \beta > -1. 
$$
Then the polynomials $Q_{j,\nu}^n$ in \eqref{eq:Q-basis} with $q_j = 
p_j^{\b_{n-2j}}$, where $\b_k=k+(d-2)/2$, form a mutually orthogonal
basis for $\CV_n^d(\Delta)$.
\end{lem}

\begin{proof}
It is easy to see that $(f,g)_\beta$ is indeed a positive definite inner
product on the space of polynomials of one variables, so that the orthogonal
polynomials with respect to $(f,g)_\beta$ exist (see Lemma \ref{lem:2.3} 
below). Using the formula
$$
 \int_{B^d} f(x) dx = \int_{0}^1 r^{d-1} \int_{S^{d-1}} f(r x') d\omega(x')dr,
$$ 
the definition of $Q_{j,\nu}^n$ and \eqref{eq:harmonics} shows immediately that 
\begin{align*}
 \langle Q_{j,\nu}^n, Q_{j',\nu'}^{n'} \rangle_\Delta  & :=  
    \delta_{\nu,\nu'} \delta_{n-2j,n'-2j'} \\    
   & \times \frac{1}{4d} \int_0^1 r^{d+2(n-2j)-1} 4^2 
     (\CJ_{\b_{n-2j}} q_j)(2r^2-1) (\CJ_{\b_{n'-2j'}} q_{j'})(2r^2-1) dr. \notag 
\end{align*}
In the nonzero case we have $\b_{n-2j} = \b_{n'-2j'}$.
Thus, a  change of variable $r \mapsto \sqrt{(1+s)/2}$ shows that  
\begin{align} \label{eq:norm}
 \langle Q_{j,\nu}^n, Q_{j',\nu'}^{n'} \rangle_\Delta  
  = \delta_{\nu,\nu'} \delta_{n-2j,n'-2j'}   \frac{1}{d 2^{\beta_{n-2j}}}
   (q_j,q_{j'})_{\b_{n-2j}},
\end{align}
which proves the stated result. 
\end{proof}

We note that $q_j^{\beta_{n-2j}}$ should be understood as one member (of 
degree $j$) in the orthogonal family $\{q_k^{\b_{n-2j}}\}$.

\begin{lem} \label{lem:2.3}
The polynomials $p_j^\beta$ defined by
$$
 p_0^\beta(s) = 1, \qquad p_j^\beta(s) = 
     (1-s) P_{j-1}^{(2,\b)}(s), \quad j \ge 1,
$$
are orthogonal with respect to the inner product $(f,g)_\beta$. 
\end{lem}

\begin{proof}
We need the following property of the Jacobi polynomials \cite[p. 71]{Sz},
\begin{equation} \label{eq:JacobiA}
 (1-s)P_{j-1}^{(2,\b)}(s) = \frac{2}{2j+\b +1}
     \left[(j+1)P_{j-1}^{(1,\b)}(s) - j P_j^{(1,\b)}(s)\right]. 
\end{equation} 
The Jacobi polynomial $P_{j-1}^{(1,\b)}$ satisfies a differential equation
$$
  (1-s^2) y'' + (-1+\b -(3+\b)s )y' +(j-1)(j+\b+1)y =0.
$$
Using these two facts, we easily deduce that 
\begin{align*}
 & \frac{2j+ \b +1}{2} \CJ_\b \left[(1-s)P_{j-1}^{(2,\b)}(s)\right] = 
    (j+1) \CJ_\b P_{j-1}^{(1,\b)}(s) - j \CJ_\b P_{j-1}^{(1,\b)}(s)\\ 
& \qquad =(j+1) \left[(-(j-1)(j+\b+1) -(\b+1)) P_{j-1}^{(1,\b)}(s) \right. \\
    & \qquad\qquad \qquad\qquad \qquad   
 \left.  - j(-j(j+\b+2)-(\b+1)) P_j^{(1,\b)}(s)\right]\\ 
& \qquad =- j(j+1) \left[(j+\b) P_{j-1}^{(1,\b)}(s) -
   (j+\b +1) P_j^{(1,\b)}(s)\right].  
\end{align*} 
We need yet another formula of Jacobi polynomials \cite[p. 782, (22.7.8)]{AS},
\begin{equation} \label{eq:Jacobi}
 (2j+\b+1) P_j^{(0,\b)}(s) = (j+\b+1)P_j^{(1,\b)}(s)- (j+\b)P_{j-1}^{(1,\b)}(s)
\end{equation}
which implies immediately that 
\begin{equation}\label{eq:JPn}
 \CJ_\b \left[(1-s) P_{j-1}^{(2,\b)}(s)\right]
 = 2j (j+1)P_j^{(0,\b)}(s). 
\end{equation}
Hence, for $j,j' \ge 1$, we conclude that 
\begin{align}\label{eq:norm2}
 (p_j^\beta, p_{j'}^\beta)_\beta = & \int_{-1}^1 
     \CJ_\b \left[(1-s) P_{j-1}^{(2,\b)}(s)\right] \CJ_\b
       \left[(1-s) P_{j'-1}^{(2,\b)}(s)\right](1+s)^\beta)ds \\
 = & 4j(j+1)j'(j'+1) \int_{-1}^1 
     P_j^{(0,\b)}(s) P_{j'}^{(0,\b)}(s) (1+s)^\b ds = 0  \notag
\end{align}
whenever $j \ne j'$. Furthermore, for $j \ge 1$, we have 
$$
 (p_0^\beta, p_j^\beta)_\beta = - 2j(j+1)(\beta+1) \int_{-1}^1 
     P_j^{(0,\b)}(s) (1+s)^\b ds = 0  
$$
since $(\CJ_\beta p_0^\beta)(s) = (\CJ_\beta 1) = - (\b +1)$. 
\end{proof}

As a consequence of the above lemmas, we have found a mutually orthogonal 
basis with respect to $\langle \cdot,\cdot \rangle_\Delta$.

\begin{thm}
A mutually orthogonal basis for $\CV_n^d(\Delta)$ is given by
\begin{align}\label{eq:basis}
\begin{split}
Q_{0,\nu}^n(x) & = Y_\nu^n(x), \quad \\
 Q_{j,\nu}^n(x) & = (1-\|x\|^2) P_{j-1}^{(2,n-2j+\frac{d-2}{2})}
     Y_\nu^{n-2j}(x), \quad 1 \le j \le \frac{n}{2},
\end{split}
\end{align}
where $\{Y_\nu^{n-2j}: 1 \le \nu \le \sigma_{n-2j}\}$ is an orthonormal basis 
of $\CH_{n-2j}^d$. Furthermore, 
\begin{align}\label{eq:Qnorm}
\langle Q_{0,\nu}^n, Q_{0,\nu}^n \rangle_\Delta = \frac{2n+d}{d}, \qquad
\langle Q_{j,\nu}^n, Q_{j,\nu}^n \rangle_\Delta =\frac{8 j^2(j+1)^2}{d(n+d/2)}.
\end{align}
\end{thm}

\begin{proof}
The fact that $Q_{j,\nu}^n \in \CV_n^d(\Delta)$ follows from Lemma 
\ref{lem:2.2} and Lemma \ref{lem:2.3}. To compute the norm of $Q_{0,\nu}^n$ 
we use the fact that 
\begin{equation} \label{eq:Euler}
   \Delta [(1-\|x\|^2)Y_\nu^{n-2j}(x)] = -2 d Y_\nu^n(x)
         -4 \langle x,\nabla\rangle Y_\nu^n = -2 (d+2n)Y_\nu^n(x)
\end{equation}
by Euler's formula on homogeneous polynomials, which shows that 
$$
\langle Q_{0,\nu}^n, Q_{0,\nu}^n \rangle_\Delta = 
  \frac{(2n+d)^2}{d} \int_0^1 r^{d-1+2n} dr \frac{1}{\omega_{d-1}}
     \int_{S^{d-1}} \left[Y_\nu^n(x)\right]^2 dx = \frac{2n+d}{d}. 
$$
Furthermore, using the equation \eqref{eq:norm} and \eqref{eq:norm2}, we have
\begin{align*}
\langle Q_{j,\nu}^n, Q_{j,\nu}^n \rangle_\Delta
 & = \frac{1}{d 2^{\beta_j}} (p_j,p_{j'})_{\b_j} =
     \frac{4j^2(j+1)^2}{d 2^{\beta_j}} 
     \int_{-1}^1\left[P_j^{(0,\beta_j)}(s)\right]^2 (1+s)^{\b_j} ds\\   
 & = \frac{8 j^2(j+1)^2}{d(\b_j+2j+1)} =\frac{8 j^2(j+1)^2}{d(n+d/2)}
\end{align*}
where we have used the well known formula for the norm of the Jacobi polynomial
(see, for example, \cite[p. 68]{Sz}).
\end{proof}

The explicit formula of the basis \eqref{eq:basis} leads to the following
interesting result, which relates $\CV_n^d(\Delta)$ to orthogonal polynomials
with respect to $W_2(x) = (1-\|x\|)^2$. 

\begin{cor}
For $n \ge 1$, 
$$
  \CV_n^d(\Delta) = \CH_n^d \oplus (1-\|x\|^2)\CV_{n-1}^d(W_2).
$$
\end{cor}

\begin{proof}
Using the basis \eqref{eq:Wmu-basis} for $\CV_{n-1}^d(W_2)$, it follows that
we actually have
\begin{equation} \label{eq:Q-P}
  Q_{j,\nu}^n(x) = (1-\|x\|^2) P_{j-1,\nu}^{n-2}(W_2;x)
\end{equation}
for $j \ge 1$, from which the stated result follows. 
\end{proof}

In the case of $d =2$, an orthonormal basis for the space $\CH_k^2$ is 
given by
$$
  Y_1^n (x,y)= \sqrt{\tfrac{1}{2}}r^n \cos n\theta \quad\hbox{and}\quad
  Y_2^n (x,y)= \sqrt{\tfrac{1}{2}} r^n \sin n \theta 
$$
in polar coordinates $x= r\cos\theta$, $y=r\sin\theta$. Hence, a mutually 
orthogonal basis for $\CV_n^2(\Delta)$ is given by 
\begin{align*}
& Q_{0,1}^n(x,y) = Y_1^n(x,y),\qquad  Q_{0,2}^n(x,y) = Y_2^n(x,y), \\
& Q_{j,1}^n(x,y)= (1-x^2-y^2)P_{j-1}^{(2,n-2j)}(2x^2+2y^2 -1)Y_1^{n-2j}(x,y), 
  \quad 1 \le j \le \tfrac{n}{2}\\ 
& Q_{j,2}^n(x,y) = (1-x^2-y^2)P_{j-1}^{(2,n-2j)}(2x^2+2y^2 -1) Y_2^{n-2j}(x,y),
  \quad 1 \le j \le \tfrac{n-1}{2},  
\end{align*}
which becomes an orthonormal basis upon dividing by the square root of the 
norm given by \eqref{eq:Qnorm}. Without normalization, this gives 
\begin{align*}
&\CV_1^2(\Delta) = \sspan \{x, y\}, \qquad 
\CV_2^2(\Delta) = \sspan \{x^2-y^2, xy,  1-x^2-y^2 \}, \\
&
\CV_3^2(\Delta)= \sspan \{ x^3 - 3 x y^2,  3y^3 - x^2 y, x (1-x^2-y^2), 
   y (1-x^2-y^2)\},
\end{align*}
for example.

\section{Expansions in Sobolev Orthogonal Polynomials}\label{Expand}
\setcounter{equation}{0}

Let $H^2(B^d)$ denote the space of functions for which 
$\langle f,f \rangle_\Delta$ is finite. This is not the $L^2$ space on
$B^d$ since the definition of $\langle \cdot,\cdot \rangle_\Delta$ require
that $f$ has second order derivatives. Nevertheless the standard Hilbert 
space theory shows that every $f\in H^2(B^d)$ can be expanded into a series
in Sobolev orthogonal polynomials. In other words, 
$$
  H^2(B^d) = \sum_{n=0}^\infty\oplus \CV_n^d(\Delta): 
        \qquad f = \sum_{n=0}^\infty \proj_n f,
$$
where $\proj_n: H^2(B^d) \mapsto \CV_n^d(\Delta)$ is the projection operator,
which can be written in terms of the orthonormal basis \eqref{eq:basis} as   
\begin{equation} \label{eq:proj}
 \proj_n f(x)=\sum_{0 \le j \le n/2} H_j^{-1}\sum_{\nu =0}^{\sigma_{n-2j}}
   \wh f_{j,\nu}^n  Q_{j,\nu}^n(x), \qquad \wh f_{j,\nu}^n = 
            \langle f,Q_{j,\nu}^n \rangle_\Delta,
\end{equation}
where $H_j = \langle Q_{j,\nu}^n, Q_{j,\nu}^n \rangle_\Delta$ are independent
of $\nu$ as shown in \eqref{eq:Qnorm}. Let $P_n^\Delta(x,y)$ 
denote the reproducing kernel of $\CV_n^d(\Delta)$. In terms of the orthonormal
basis \eqref{eq:basis} in the previous section, the reproducing kernel can be
written as 
$$
  P_n^\Delta (x,y) = 
     \sum_{0\le j \le n/2} H_j^{-1} 
       \sum_\nu Q_{j,\nu}^n(x) Q_{j,\nu}^n(y).  
$$
The projection operator can be written as an integral operator with 
$P_n^\Delta$ as its kernel, which means that 
\begin{align*}
 \proj_n f(x) & = \langle f,P_n^\Delta(x,\cdot) \rangle_\Delta \\
 & = \frac{1}{4 d^2 \vol(B^d)}
  \int_{B^d}\Delta[(1-\|y\|^2) f(y)] \Delta[(1-\|y\|) P_n^\Delta(x,y)] dy  
\end{align*}
where $\Delta$ is applied on $y$ variable. 

It turns out that the orthogonal expansion can be computed without involving
derivatives of $f$. 

\begin{prop}
For $j \ge 1$, let $\b_j = n-2j+(d-2)/2$; then
\begin{align} \label{eq:Fourier}
 \wh f_{j,\nu}^n = &\frac{8 j (j+1)}{d^2 \vol(B^d)} \left[
 (\b_j+j)(\b_j+j+1) \int_{B^d} f(x) Q_{j,\nu}^n(x)dx \right. \\
 &  \hspace{1in} 
 \left.  - \frac{1}{2} \int_{S^{d-1}} f(y')Y_\nu^{n-2j}(y') d\omega(y')\right];
\notag
\end{align} 
furthermore, for $j =0$, 
$$
 \wh f_{0,\nu}^n = \frac{d+2n}{d} 
   \frac{1}{\omega_d}\int_{S^{d-1}} Y_\nu^{n-2j}(y')f(y') d\omega(y').
$$
\end{prop}

\begin{proof}
By \eqref{eq:Wmu-basis}, $P_{j,\nu}^n(W_0;x) = P_j^{(0,\beta_j)}(2\|x\|^2-1) 
Y_\nu^{n-2j}(x)$. Let $j \ge 1$. By Lemma \ref{lem:2.1} and \eqref{eq:JPn}, 
$$
 \Delta \left[(1-\|x\|^2 Q_{j,\nu}^n(x)\right] = 
    8 j(j+1) P_{j,\nu}^n(W_0;x).
$$
Applying Green's identity 
$$
  \int_{B^d} ( u \Delta v - v \Delta u) dx = 
      \int_{S^{d-1}} \left( \frac{\partial v}{\partial n} u -  
          \frac{\partial u}{\partial n} v \right) d\omega  
$$
with $v(x) = (1-\|x\|^2)f(x)$ and $u = Q_{j\b}^n$ shows then 
\begin{align} \label{eq:3.x}
 \wh f_{j,\nu}^n & = \frac{8 j (j+1)}{4 d^2 \vol(B^d)} \int_{B^d} 
      \Delta \left[ (1-\|x\|^2)f(x) \right] P_{j,\nu}^n(W_0;x) dx \\
   & = \frac{2 j (j+1)}{ d^2 \vol(B^d)} \left [ 
        \int_{B^d} (1-\|x\|^2) f(x) \Delta P_{j,\nu}^n(W_0;x)dx \right .
    \notag\\
   & \qquad\qquad  \qquad\qquad  \qquad  
       \left.  - 2 \int_{S^{d-1}} Y_\nu^{n-2j}(x') f(x') d\omega \right]\notag
\end{align}
where we have used the fact that $P_j^{(0,\b)}(1) =1$. Let $\partial 
P^{(0,\beta)}$ denote the derivative of $P^{(0,\beta)}$. Using 
\eqref{eq:Delta} and \eqref{eq:Eigen} it is easy to see that 
\begin{align*}
\Delta[P_{j,\b}^n(W_0;x)] = &
     8\left[2r^2 \partial^2 P_j^{(0,\beta_j)}(2r^2-1) \right. \\
     & \qquad 
  \left. +(n-2j+d/2)\partial P_j^{(0,\beta_j)}(2r^2-1)\right]Y_\nu^{n-2j}(x).
\end{align*}
Let us denote the expression in the square bracket by $M_j$. The Jacobi 
polynomial $P_j^{(0,\beta)}(s)$ satisfies the differential equation
$$
  (1-s^2) y'' - (-\b + (\b +2)s ) y' + j(j+\b+1) y =0.  
$$
Hence, changing variable $2r^2 -1 \mapsto s$, we conclude that  
$$
2(1-r^2) M_j =  - j (j+\beta_j+1) P_j^{(0,\b_j)}(s) +\frac{1}{2}
       (j+\b+1)(1+s) P_{j-1}^{(1,\b_j+1)}(s).
$$
On the other hand, using \eqref{eq:Jacobi}, \eqref{eq:JacobiA} and the fact 
that \cite[p. 782]{AS}  
\begin{align*}
 (2j+\b +1) (1+s) P_{j-1}^{(1,\b+1)}(s) = 2(j+\b)P_j^{(1,\b)}(s) + 2 j
    P_{j-1}^{(1,\b)}(s)
\end{align*}
we conclude that  
\begin{align*} 
2(1-r^2) M_j & =  \frac{(\b_j+j+1)(\b_j+j)}{2j+\b_j+1}
   \left[-j P_{j}^{(1,\b_j)}(s)  + (j+1) P_{j-1}^{(1,\b_j)}(s)\right]\\  
& = \frac{1}{2} (\b_j+j+1)(\b_j+j)(1-s) P_{j-1}^{(2,\b_j)}(s) \\
& =  (\b_j+j+1)(\b_j+j)(1-r^2) P_{j-1}^{(2,\b_j)}(2r^2-1).  
\end{align*}
Consequently, we have proved that 
\begin{align*} 
(1-\|x\|^2)\Delta[P_{j,\nu}^n(W_0;x)] & = 4 (\b_j+j+1)(\b_j+j)(1-r^2)
    P_{j-1}^{(2,\b_j)}(2r^2-1) Y_\nu^{n-2j}(x) \\
 & =  4(\b_j+j+1)(\b_j+j) Q_{j,\nu}^n(x)
\end{align*} 
which leads to the stated result for $j \ge 1$ by \eqref{eq:3.x}. The proof
of $j =0$ is similar but easier, in which we need to use \eqref{eq:Euler}. 
\end{proof} 

Let us denote by $P_n(W_\mu; x,y)$ the reproducing kernel of $\CV_n^d(W_\mu)$,
which can be written as 
$$
 P_n(W_\mu; x,y) = \sum_{|\alpha| = n} A_{\alpha,\mu}^{-1} 
   P_\alpha(W_\mu;x) P_\alpha(W_\mu;y), 
$$
where $A_{\alpha,\mu} = c_\mu \int_{B^d}[P_\alpha(W_\mu;y)]^2 W_\mu(y)dy$ in
which $c_\mu$ is the normalization of $W_\mu$.
Let us also denote by $C_n^\lambda(t)$ the Gegenbauer polynomial of degree 
$n$, and by $x \cdot y$ the usual dot product of $x, y \in \RR^d$. 

\begin{cor}
For $f \in H^2(B^d)$ and $x \in B^d$,
\begin{align*}
 \proj_n f(x)& = Y_n f(x) + 
  (1-\|x\|^2)  \frac{4}{\binom{d}{2}\vol(B^d)}
   \int_{B^d} f(y) P_{n-2}(W_2;x,y) (1-\|y\|^2) dy \\
    & -\frac{(n+d/2)}{4} (1-\|x\|^2)
    \sum_{1 \le j \le n/2} \frac{P_{j-1}^{(2,n-2j+\frac{d-2}{2})}(2\|x\|^2-1)}
       {P_{j-1}^{(2,n-2j+\frac{d-2}{2})}(1)} Y_{n-2j}f(x),
\end{align*} 
where with $x' = x/\|x\| \in S^{d-1}$, 
$$
  Y_m f(x) =   \|x\|^m \int_{S^{d-1}} f(y') 
      \frac{m+(d-2)/2}{(d-2)/2} C_m^{\frac{d-2}{2}}(x \cdot y') d\omega(y').
$$ 
\end{cor} 

\begin{proof}
The values of $H_j = \langle Q_{j,\nu}^n, Q_{j,\nu}^n \rangle_\Delta$ are 
given in \eqref{eq:Qnorm}. It follows immediately that 
$$
\sum_{\nu=1}^{\sigma_n} H_0^{-1} \wh f_{0,\nu}^n Q_{0,\nu}^n(x)  
= \frac{1}{\omega_{n-1}} \int_{S^{d-1}} f(y') \sum_{\nu=1}^{\sigma_n} 
     Y_\nu^n(y')Y_\nu^n(x) d\omega(y') = Y_n f(x), 
$$
where the last step follows from the summation formula of spherical 
harmonics,
$$
 \sum_{\nu=1}^{\sigma_n} Y_\nu^n(x)Y_\nu^n(y)=
   \|x\|^n \sum_{\nu=1}^{\sigma_n} Y_\nu^n(x')Y_\nu^n(y)
   = \|x\|^n \frac{n+(d-2)/2}{(d-2)/2}
     C_n^{\frac{d-2}{2}} (x' \cdot y)
$$
for $x', y \in S^{d-1}$. Furthermore, setting $f = Q_{j,\nu}^n$ with $j \ge 1$
in \eqref{eq:Fourier} also shows
$$
 H_j = \frac{8j(j+1)}{d^2 \vol(B^d)} (\b_j+j)(\b_j+j+1)
    \int_{B^d} [ Q_{j,\nu}^n(x)]^2 dx. 
$$
Hence, it follows from \eqref{eq:Fourier} and \eqref{eq:Qnorm} that 
\begin{align}\label{eq:3.3}
H_j^{-1} \wh f_{j,\nu}^n 
=  \frac{\int_{B^d} f(y) Q_{j,\nu}^n(y)dy}{\int_{B^d} [Q_{j,\nu}^n(y)]^2dy} 
    - \frac{n+d/2}{2 j(j+1)} 
      \frac{1}{\omega_{d-1}} \int_{S^{d-1}} f(y')Y_\nu^{n-2j}(y') d\omega(y').
\end{align}
The relation \eqref{eq:Q-P} readily shows that 
\begin{align}\label{eq:3.4}
  \int_{B^d} [Q_{j,\nu}^n(y)]^2dy = \frac{1}{4} \binom{d}{2} \vol(B^d)
     c_2 \int_{B^d} [P_{j,\nu}^n(W_2;y)]^2 
      W_2(y) dy.
\end{align}
We multiply \eqref{eq:3.3} by $Q_{j,\nu}^n(x)$ and sum over $\nu$ 
and $j$. Using \eqref{eq:3.4} and the fact that 
$P_{j-1}^{(2,n-2j+\frac{d-2}{2})}(1)= j(j+1)/2$, the stated result follows 
from \eqref{eq:Fourier} and \eqref{eq:proj}.
\end{proof}

It follows from this corollary that the orthogonal expansion of $f$ with 
respect to $\langle \cdot,\cdot \rangle_\Delta$ coincides with the spherical 
harmonic expansion of $f$ when restricted on $S^{d-1}$.

\bigskip\noindent
{\it Acknowledgment.} The author thanks Professor Ken Atkinson for drawing
his attention to this problem.

\end{document}